\documentclass[12pt]{article}
\usepackage{latexsym, amsmath, amsfonts, amscd, amssymb, verbatim}
\usepackage{enumerate}
\def\tto{\;{\lower 1pt \hbox{$\rightarrow$}}\kern -10pt
\hbox{\raise 2pt \hbox{$\rightarrow$}}\;}

\def\emp{\emptyset}

\newcommand{\qedsymbol}{$\blacksquare$}
\begin{document}
\pagestyle{myheadings}

\newtheorem{Theorem}{Theorem}[section]
\newtheorem{Proposition}[Theorem]{Proposition}
\newtheorem{Remark}[Theorem]{Remark}
\newtheorem{Lemma}[Theorem]{Lemma}
\newtheorem{Corollary}[Theorem]{Corollary}
\newtheorem{Definition}[Theorem]{Definition}
\newtheorem{Example}[Theorem]{Example}
\renewcommand{\theequation}{\thesection.\arabic{equation}}
\normalsize

\setcounter{equation}{0}

\title{\bf  Constrained Fermat-Torricelli-Weber Problem in real Hilbert Spaces}

\author{Son Dang Nguyen}\maketitle

\medskip
\begin{quote}
\noindent {{\bf Abstract} The Fermat-Weber location problem requires finding a point in $\mathbb{R}^n$ that minimizes the sum of weighted Euclidean distances  to $m$ given points. An iterative solution method for this problem was first introduced by E. Weiszfeld in 1937. Global convergence of Weiszfeld's algorithm was proven by W. Kuhn in 1973. This paper studies Fermat-Weber location problems with closed convex constraint sets in real Hilbert spaces. In section two, we show that existence and uniqueness of solutions of the problems. Moreover, the solution is stable with respect to the perturbation of the $m$ anchor points. In the section three, we extend Weiszfeld's algorithm by adding a projection on the constraint set. The convergence of the sequence generated by the method to the optimal solution of the problem is proved. .}\  


\medskip
\noindent {{\bf Key Words} \text{ . } Fermat-Weber problem in Hilbert spaces \text{ . } Closed convex constraint \text{ . } Extension of solution stability \text{ . } Convergence of Weiszfeld's method \text{ . } }\ 
\end{quote}

\section{Introduction and Preliminary}
\markboth{\centerline{\sc Introduction}}{\centerline{\sc N.D.Son}} \setcounter{equation}{0}
Fermat's location problem from the seventeenth century is stated as follows: Given three anchor points in a plane, find a fourth point such that the sum of its distances to the three given anchor points is as small as possible. The Italian physicist and mathematician E. Torricelli found a method to construct the unique solution point, which later was called {\it the Fermat-Toricelli point}. At the beginning of the twentieth century, Weber (see \cite{Weber}), studied location problems with weights and with more than three anchor points. The new problem was consequently called {\it the generalized Fermat-Weber location problem}. It also bears other names like {\it the Fermat problem, the Weber problem, the Fermat-Toricelli problem, the Steiner problem, etc.}. The generalized Fermat-Weber location problem is the following non-smooth, convex optimization problem
\begin{equation}\label{11}
\min \left\{ {f(x) = \sum\limits_{i = 1}^m {{w_i}\left\| {x - {a_i}} \right\|} \left| x \in \mathbb{R}^n\right.} \right\},
\end{equation}
where $m$ points $a_1,a_2,...,a_m$ are given in $\mathbb{R}^n$, called the anchors and $w_1,w_2,...,w_m$ are positive numbers, called weights. It is well known that if the anchor points are not collinear, i.e there does not exist any straight line containing all the points $a_1,a_2,...,a_m$, then the  objective function of \eqref{11} is strictly convex and coercive. The function $f(x)$ tends to $+\infty$ as $\left\| x \right\| \to  + \infty $ . (If $a_1,a_2,...,a_m$ are collinear then at least one of the points $a_1,a_2,...,a_m$  is optimal and it can be found in linear time, see \cite{Blum_Floyd_Pratt_Rivest_Tarjan}). To solve \eqref{11}, there were several schemes \cite{J.A. Chatelon D.W.Hearn and T.J. Lowe, M.L. Overton} and presented by Weiszfeld in \cite{Weiszfeld1,Weiszfeld2} was the one of the most popular methods. The results of Weiszfeld were rediscovered several years later independently  by Miehle  \cite{Miehle},  Kuhn and Kuenne  \cite{Kuhn and Kuenne}, and Cooper \cite{L. Cooper}. Weiszfeld's algorithm is based on the next mapping $T:\mathbb{R}^n \to \mathbb{R}^n$,
\begin{equation} \label{1.1}
T(x) = \left\{ \begin{array}{l}
\begin{array}{*{20}{c}}
{\dfrac{1}{{\sum\nolimits_{i = 1}^m {\dfrac{{{w_i}}}{{\left\| {x - {a_i}} \right\|}}} }}\sum\limits_{i = 1}^m {\dfrac{{{w_i}{a_i}}}{{\left\| {x - {a_i}} \right\|}}} }&{\text{if } x \ne a_1,a_2,...,a_m}
\end{array}\\
\\
\begin{array}{*{20}{c}}
{{a_j}}&{\text{if } x = {a_j} \text{ for some $j$}}
\end{array}
\end{array} \right..
\end{equation}
Weiszfeld's algorithm is defined by the iterative scheme: 
\begin{equation} \label{1.2}
{x_{k + 1}} = T({x_k}), \text{} k=0,1,2,... \text{ . }\
\end{equation}

Convergence of the above algorithm are discussed in \cite{Katz,Kuhn}. In 1973,
Kuhn \cite{Kuhn} claimed that $\{x_k\}$ converges to the unique solution for
all but a denumerable number of starting points $x_0$. However, Chandrasekaran and Tamir \cite{R. Chandrasekaran and A. Tamir} detected a flaw in the Kuhn's statement and showed that the system $T(x)=a_i$ may have a continuum set of solutions even when the points $a_1,...,a_m$ are not
collinear. Brimberg \cite{Brimberg} proved the conjecture of
Chandrasekaran and Tamir, but in \cite{Canovas and Canavate} Canovas et
al. found counterexamples to the proof in Brimberg's paper. Eight years later
Brimberg modified the previous proof and solved the proof and solved the conjecture of Chandrasekaran and Tamir. Finally, Canavate
\cite{Canavate} claimed that the Weiszfeld's algorithm converges
for all points but a set of measure zero.

In this work, we will consider closed convex constrained location problem,
\[\min \left\{ {\left. {f(x) = \sum\limits_{i = 1}^m {{w_i}\left\| {x - {a_i}} \right\|} } \right|x \in C} \right\},\begin{array}{*{20}{c}}
{}&{}&{}&{}
\end{array}\left( {FTW} \right)\]

We mention some concepts about normal cone of convex sets and subgradient of convex function. Normal cone of a closed convex set, notes $N(x,C )$, is defined 
\[N(x,C) = \left\{ {\begin{array}{*{20}{l}}
	{\begin{array}{*{20}{c}}
		{\{ {y^*} \in H\left| {\left\langle {\text{ }{y^*},z - x} \right\rangle  \le 0,{\text{\rm{ for all }}}z \in C} \right.\} }&{{\text{\rm{ if }}}x \in C}
		\end{array}}\\
	{\begin{array}{*{20}{c}}
		\emptyset &{{\text{\rm{ if }}}x \notin C}
		\end{array}}
	\end{array}} \right..\]

Let $f: H \to {\bar {\mathbb{R}}}$ be a convex function and let $\bar{x} \in dom  \text{ \it f}$. A subgradient of $f$ at $\bar{x}$ is noted $\partial f(\bar x)$ if 
\[\partial f(\bar x) = \left\{ {\left. {{y^*} \in H} \right|\text{} \left\langle {{{y^*}},x - \bar x} \right\rangle  \le f(x) - f(\bar x)} \text{ for all x} \in H\right\} .\]
It is easy that $\partial {I_C}(x) = N(x,C)$ with for all $x \in C$.
where $C$ is a closed convex subset of a real Hilbert space $H$, the anchor points ${a_1},{a_2}, \ldots ,{a_m}$ are in  $H$, and $w_i,$ $ i=1,2,...,m$ are positive weights. In this paper, we always assume that the anchor points $a_1,a_2,...,a_m$ are collinear and note that the points ${a_1},{a_2}, \ldots ,{a_m}$ can belong to $C$ or not. If all the points is in $C$, then we have a non-constrained problem (FTW).

 The distance function associated with a nonempty, closed convex set $C \subset H$ is,
\[d(x;C) = \inf \{ \left. {\left\| {x - y} \right\|} \right|y \in C\} .\]

\begin{Remark} \label{remark1.1}\rm
	The distance function  has the following properties:
	\begin{enumerate}[(a)]
		\item It is a Lipschitz continuous convex function.
		\item Let ${\Pi _C }(x) = \{ y \in C \left| {\left\| {x - y} \right\| = d(x,C )} \right.\} $, then ${\Pi _C }(x) $ has a unique point. The point $y = {\Pi _C }(x)$ is called the projection of $x$ on $C$. Note that ${\Pi _C }(x)=x $, for all $x \in C$.
		\item See [D. Kinderlehrer and G. Stampacchia, Theorem 2.3 \cite{D.K and G. S}], for every $x \in C$, then $y = {\Pi _C }(x)$ if and only if
		\[\begin{array}{*{20}{c}}
		{\begin{array}{*{20}{c}}
			{y \in C:}&{\left\langle {x - y,z - y} \right\rangle  \le 0}
			\end{array}}&{\text{for all $z \in C$}}
		\end{array}\]
		or
		\[ x - {\Pi _C }(x)\in N({\Pi _C }(x),C ).\]	
		
		\item The mapping projection: ${\Pi _C }(.):H \to C$ is continuous and nonexpansive, or
		\[\left\| {{\Pi _C}(x) - {\Pi _C}({x'})} \right\| \le \left\| {x - {x'}} \right\| \text{ for }x,{x'}\in H.\]
	\end{enumerate}
\end{Remark}
\section{Analysis}	
\subsection{The existence and uniqueness solution}
The problem (FTW) is a convex optimization with closed convex constrained, so it is possible to rewrite with non-constrained,
\[\min \left\{ {f(x) = \sum\limits_{i = 1}^m {{w_i}\left\| {x - {a_i}} \right\|}  + {I_C}(x)\left| {x \in H} \right.} \right\},\begin{array}{*{20}{c}}
{}&{}&{}&{}
\end{array}\left( {FTW} \right)\]
where ${I_C }(x)$ is indicator function, is defined
\[{I_C }(x) = \left\{ \begin{array}{l}
\begin{array}{*{20}{c}}
0&{\text{if } x \in C }
\end{array}\\
\begin{array}{*{20}{c}}
{ + \infty }&{\text{if } x \notin C }
\end{array}
\end{array} \right..\]

The objective function of the problem (FTW) is a continuous convex function in norm-topology of the real Hilbert space H, so it is lower semi-continuous function in weak-topology. It also satisfies coercive condition, hence the solution of the problem (FTW) exists. Although it is well-known fact that the function is strictly convex in $\mathbb{R}^n$ if the points $a_1,a_2,...,a_m$ are \textbf{not collinear} \cite{Kuhn,Kuhn and Kuenne,Weiszfeld1,Weiszfeld2} i.e. all $a_i,i=1,...,m$ do not lie on a certain straight line, we will still prove it again similarly in the real Hilbert space $H$.
\begin{Lemma}
If the points $a_1,a_2,...,a_m$ are not collinear, the objective function of (FTW) is strictly convex in a real Hilbert space, and hence so is it in closed convex constraint sets. 
\end{Lemma}
{\bf Proof}. Since each function ${f_i}(x): = \left\| {x - {a_i}} \right\|$ as $i=1,2,...,m$ is obviously convex, their sum $f = \sum\limits_{i = 1}^m {{f_i}} $ as well, i.e., for any $x,y \in H$ and $\lambda \in (0,1)$ we have
\begin{equation} \label{4.2}
f\left( {\lambda x + (1 - \lambda )y} \right) \le \lambda f(x) + (1 - \lambda )f(y).
\end{equation}
 Supposing by contradiction that $f$ is not strictly convex, find $\bar x,\bar y \in H$ with $\bar x \ne \bar y$ and $\lambda \in (0,1)$ such that \eqref{4.2} holds as equality. It follows that
 \[{f_i}\left( {\lambda x + (1 - \lambda )y} \right) = \lambda {f_i}(x) + (1 - \lambda ){f_i}(y) \text{ for all } i=1,2,...m,\]
 which ensures therefore that
 \[\left\| {\lambda (\bar x - {a_i}) + (1 - \lambda )(\bar y - {a_i})} \right\| = \left\| {\lambda (\bar x - {a_i})} \right\| + \left\| {(1 - \lambda )(\bar y - {a_i})} \right\|,\begin{array}{*{20}{c}}
 {}&{i = 1,2,...,m.}
 \end{array}\]
 If $\bar x \ne {a_i}$ and $\bar y \ne {a_i}$, then there exists $t_i >0$ such that ${t_i}\lambda (\bar x - {a_i}) = (1 - \lambda )(\bar y - {a_i}),$ and hence
 \[\bar x - {a_i} = {\gamma _i}(\bar y - {a_i})\text{ with }{\gamma _i} = \frac{{1 - \lambda }}{{{t_i}\lambda }}.\]
 Since $\bar x \ne \bar y$, we have ${\gamma _i} \ne 1.$ Thus
 \[{a_i} = \frac{1}{{1 - {\gamma _i}}}\bar x - \frac{{{\gamma _i}}}{{1 - {\gamma _i}}}\bar y \in L(\bar x,\bar y),\]
 where $L(\bar x,\bar y)$ signifies the line connecting $\bar x$ and $\bar y$. Both cases where $\bar x=a_i$ and $\bar y = a_i$ give us $a_i \in L(\bar x, \bar y)$. Hence $a_i \in L(\bar x, \bar y)$ for $i=1,2,...,m$, a contradiction.  From the strictly convex property of the function, (FTW) has unique solution. \qedsymbol
 
Basing on basic subdifferential calculus of convex analysis in real infinite Hilbert space in the book \cite{Bauschke, Heinz H. Combettes, Patrick L.} we have easily necessary and sufficient conditions theorem such that ${{\bar x} }$ is the optimal solution of the problem (FTW) in a real Hilbert space $H$.
\begin{Theorem} \label{theorem2.1}
	Note $A = {\rm{\{ }}{{\rm{a}}_1},{a_2}, \ldots ,{a_m}\}$ and suppose the vertices $a_1,a_2,...,a_3$ are not collinear. Then ${{\bar x} }$ is the unique optimal solution of the problem (FTW) if only if 
	$$0 \in \partial f({{\bar x}}),$$
	
	Only two possible cases follows:
	\begin{enumerate}[i)]
		\item ${\bar x} \notin A$ if and only if 
		\begin{equation} \label{2.5}
		{\Pi _C} \circ T({\bar x})={\bar x}.
		\end{equation}
		or it is equivalent
		\begin{equation} \label{2.6}
		\left\langle {T(\bar x) - \bar x,x - \bar x} \right\rangle  \le 0 \text{  for every $x \in C$}
			\end{equation}
		\item ${\bar x} = a_j$, for some $j \in {1,2,...,m},$ if and only if there exists $u \in {{\bar B}_{({0},1)}}$ such that
		\begin{equation} 
	\left\langle { - \sum\limits_{i = 1,i \ne j}^m {\frac{{{a_j} - {a_i}}}{{\left\| {{a_j} - {a_i}} \right\|}}} -u ,x - {a_j}} \right\rangle  \le 0 \begin{array}{*{20}{c}}
	{}&{\text{for all $x \in C$}}.
	\end{array}
		\end{equation}
In the case, additionally $C$ is a cone or $int(C) \ne \emp$.  Then ${\bar x} = a_j$ if and only if
\[\left\| {\sum\limits_{i = 1,i \ne j}^m {\frac{{{a_j} - {a_i}}}{{\left\| {{a_j} - {a_i}} \right\|}}} } \right\| \le 1\]
	\end{enumerate} 
\end{Theorem}
{\bf Proof.} 
To prove i), assume that ${{\bar x} }$ is the unique optimal solution of the problem (FTW). We begin by writing the optimality condition $0 \in \partial f({{\bar x} })$ with noticing ${\bar x} \in C$ and ${\bar x} \notin A = {\rm{\{ }}{{\rm{a}}_1},{a_2}, \ldots ,{a_m}\} $:
\[0 \in \partial f({\bar x}) \Leftrightarrow 0 \in \sum\limits_{i = 1}^m {{w_i}} \dfrac{{{\bar x} - {a_i}}}{{\left\| {{\bar x}- {a_i}} \right\|}} + N({\bar x},C),\]

Note that ${\sum\nolimits_{i = 1}^m {\dfrac{{{w_i}}}{{\left\| {{\bar x} - {a_i}} \right\|}}} }>0$ and $\lambda N({{\bar x} },C) = N({{\bar x} },C)$ with $\lambda >0$. So
\[\dfrac{1}{{\sum\nolimits_{i = 1}^m {\dfrac{{{w_i}}}{{\left\| {{\bar x} - {a_i}} \right\|}}} }}\sum\limits_{i = 1}^m {\dfrac{{{w_i}{a_i}}}{{\left\| {{\bar x} - {a_i}} \right\|}} \in {\bar x} + } N({\bar x},C ),\]
or
\begin{equation} \label{2.1}
T({\bar x}) \in {\bar x} + N({\bar x},C ),
\end{equation}
where the mapping $T: H \to H$ is defined as \eqref{1.1}. From part c) of \ref{remark1.1} , the formula \eqref{2.1} equivalents i).

For the proof ii), if $\bar x=a_j$, by basic calculus in convex analysis we have
 \[0 \in \partial f({\bar x}) = \sum\limits_{i = 1,j \ne i}^m {{w_i}} \frac{{{a_j} - {a_i}}}{{\left\| {{a_j} - {a_i}} \right\|}} + {{\bar B}_{({0},1)}} + N({a_j},C),\]
 where ${{\bar B}_{({0},1)}}$ is unit circle in real $H$ space. It can be rewrited that there exists $u \in {{\bar B}_{({0},1)}}$ such that
\begin{equation}\label{2.9}
  - \sum\limits_{i = 1,j \ne i}^m {{w_i}} \frac{{{a_j} - {a_i}}}{{\left\| {{a_j} - {a_i}} \right\|}} - u \in N({a_j},C),
 \end{equation}
 it is equivalent
 \[\left\langle { - \sum\limits_{i = 1,i \ne j}^m {\frac{{{a_j} - {a_i}}}{{\left\| {{a_j} - {a_i}} \right\|}}}  - u,x - {a_j}} \right\rangle  \le 0\begin{array}{*{20}{c}}
 {}&{\text{for all $x \in C$}}.
 \end{array}\]
Finally, when $C$ is a cone or $int(C) \ne \emp$, then ${ - \sum\limits_{i = 1,i \ne j}^m {\frac{{{a_j} - {a_i}}}{{\left\| {{a_j} - {a_i}} \right\|}}}  - u}=0$ and thus the proof is completely proved. \qedsymbol
\subsection{Stability of the solution of the problem (FTW)}
An interesting question is discussed that the unique solution of the problem (FTW) is stable? The stability means that if we perturb the vectors $a_1,a_2,...,a_m$ then the new solutions of the problem (FTW) is near to the initial solution. The answer of the above-question is true. We will prove this conclusion for closed-convex-constrained Fermat-Weber problem. This proof for the non-constrained Fermat-Weber problem is similar. The problem (FTW) depending the point $a$ can be observe as follows
\[\min \left\{ {\sum\limits_{i = 1}^m {w_i}{\left\| {x - {a_i}} \right\|}  + {I_C}(x)\left| {a \in } \right.\Delta } \right\}\begin{array}{*{20}{c}}
{}&{(FTW_a)} 
\end{array}\]
where  the set $\Delta  := \left\{ {a = ({a_1},{a_2},...,{a_m}) \in {H^m}\left| {\text{all the points $a_1,a_2,...,a_m \in H$ are not colinear }} \right.} \right\}.$ Define the function $g:H \times {H^m} \to \mathbb{R}$,
\[g(x,a) = \sum\limits_{i = 1}^m {w_i}{\left\| {x - {a_i}} \right\|} + {I_C}(x)\]
where the Euclidean norm is in $H$ as (FTW). Observe that 
\[\begin{array}{*{20}{l}}
{{H^{m + 1}}\mathop  \to \limits^A {H^{m + 1}}\mathop  \to \limits^h {\mathbb{R}}}\\
{g = h \circ A}\\
{A({x_1},...,{x_{m + 1}}) = ({x_1},{x_1} - {x_2},{x_1} - {x_3},...,,{x_1} - {x_{m + 1}})}\\
{h({y_1},...,{y_m}) = \sum\limits_{i = 2}^{m + 1} {{w_{i-1}}} \left\| {{y_i}} \right\| + {I_C}({y_1})}
\end{array}\]

Remark that $A$ is a \textit{surjective continuous linear operator} and $h$ is a \textit{continuous convex function}. Therefore $g(.,.)$ is a \textit{continuous convex function }on its domain. It is clear that  $\text{\rm dom }g = C \times H^m$ but we wil figure out the properties of $g(.)$ on $C \times \Delta \subset  C \times H^m$.


 The function $m:\Delta \subset H^n \to \mathbb{R}$
 \begin{equation} \label{optimalvaluefunction}
 m(a): = \inf \left\{ {g(x,a) = \sum\limits_{i = 1}^m {w_i}{\left\| {x - {a_i}} \right\|}  + {I_C}(x)\left| {x \in H} \right.} \right\}
  \end{equation}
 is called the {\it optimal value function}. Since $g$ is convex, so \textbf{$m$ is convex} as well. The set
 \begin{equation}\label{solutionmapping}
 M(a): = \left\{ {x \in H\left| {g(x,a) = } \right.m(a)} \right\}
 \end{equation}
 is called the \textit{solution set} of $(FTW_a)$. From the existence and uniqueness solution of (FTW) with for all $a \in \Delta$, we see that the solution set $M(a)$ only consists of a point. Therefore we can observe one like the \textit{solution mapping} $M: \Delta \subset H^n\to H$. The following remark is important to prove the continuous property of te solution mapping $M(.)$.
 
 \begin{Remark} \label{Remark2.4}
 	 $M(a) \in {\Pi _C} (co\{ {a_1},...,{a_m}\}) $ with $a=(a_1,...,a_m)$ since $T(x) \in co\{ {a_1},...,{a_m}\} $ for all $x \in H$. 
 \end{Remark}

 Some basis definitions and theorems of functions analysis is necessarily showed,(Rudin, \cite{Rudin}).
 \begin{Definition}
 	A subset set $E$ of a Banach space $X$  is called totally bounded if $E$ lies in the union of finitely many open balls of radius $\epsilon$, for every $\epsilon > 0$.
 \end{Definition}
 
 \begin{Theorem}
 	A subset $E$ of a Banach space $X$ is totally bounded, then the closed convex hull of $E$, written $\overline {co} (E)$, is the compact convex set.
 \end{Theorem}
 
 We have the following main theorem about the stability of the solution of the problem (FTW).
\begin{Theorem}
	The solution mapping $M:H^n \to H$ in \eqref{solutionmapping} is continuous on $\Delta$, therefore so the optimal value function $m: H^n \to \mathbb{R}$ in 	\eqref{optimalvaluefunction} is locally Lipschitz convex on $\Delta$.  
	
\end{Theorem}
{\bf Proof.} Take any $a=\{a_1,a_2,...,a_m\} \in \Delta$ and ${a^k} = \{ a_1^k,a_2^k,...,a_m^k\} $ with $k \ge 0$ is the sequence which converges to $a$, then the sequence ${\{ M({a^k})\} _{k \ge 0}}$ converges to $M(a)$. Indeed, we observe the subset $M$ of $H$, being defined by
\[\Omega : = \left( {\bigcup\limits_{i = 1}^m {{{\{ a_i^k\} }_{k \ge 0}}} } \right)\bigcup {\left\{ {{a_1},{a_2},...,{a_m}} \right\}}.\]
The set $\Omega$ is totally bounded since the sequence ${\{ a_i^k\} _{k \ge 0}}$ converges to $a_i$ with $i=1,2,...,m$ (converging in the norm of $H$). It implies that $\overline {co} (\Omega)$ is the compact set. From \ref{Remark2.4}, we see that the sequences ${\{ M({a^k})\} _{k \ge 0}} \subset {\Pi _C}(\overline  {co} (\Omega))$, so there exists the convergence subsequence ${\{ M({a^{{k_l}}})\} _{l \ge 0}}$ (convergence in the norm of the real Hilbert space $H$), supposing that it converges to some point $y$. We will prove that $y=M(a)$. Indeed, for every $k \ge 0$ we always get,
\begin{equation} \label{2.11}
m({a_k}) = \sum\limits_{i = 1}^m {{w_i}\left\| {M({a^k}) - a_i^k} \right\|}  \le \sum\limits_{i = 1}^m {{w_i}\left\| {M(a) - a_i^k} \right\|}  = g(M(a),{a_k}),
\end{equation}
We observe the formula \eqref{2.11} with regard to the subsequence ${\{ M({a^{{k_l}}})\} _{l \ge 0}}$,
\begin{equation} \label{2.12}
m({a_{{k_l}}}) = \sum\limits_{i = 1}^m {{w_i}\left\| {M({a^{{k_l}}}) - a_i^{{k_l}}} \right\|}  \le \sum\limits_{i = 1}^m {{w_i}\left\| {M(a) - a_i^{{k_l}}} \right\|}  = g(M(a),{a_{{k_l}}}),
\end{equation}
According to that the sequences ${\{ a_i^k\} _{k \ge 0}}$  converges to $a_i$  with $i=1,2,...,m$ and ${\{ M({a^{{k_l}}})\} _{l \ge 0}}$ converges to $y$, respectively, and the continuity of the norm, so \eqref{2.12} becomes
\[g(y,a) = \sum\limits_{i = 1}^m {{w_i}\left\| {y - {a_i}} \right\|}  \le \sum\limits_{i = 1}^m {{w_i}\left\| {{M(a)} - {a_i}} \right\|}  = g({M(a) },a).\]
From the unique solution of the problem (FTW), $y$ exactly equals $M(a)$. It is clear that the mapping $a \to (M(a),a)$ is continuous. Thus the continuity of the optimal value function $m(.)$ is easily found. Notice that $m()$ is a continuous convex on $\Delta$, thus it is locally Lipschitz that may be seen in \cite{Bauschke Heinz H. Combettes Patrick L.} \qedsymbol

The following theorem is a formula of convex subdifferential of optimal value function $m()$.

\begin{Theorem}
Subdifferential of optimal value function $m()$ can be exhibited as follows:
\[\begin{array}{ll}
\partial m(a) &= \left\{ {{a^*} \in {H^m}\left| {(0,{a^*}) \in \partial g(M(a),a)} \right.} \right\}\\
\\
&= \left\{ \begin{array}{l}
\left[ {\begin{array}{*{20}{c}}
	{ - \frac{{{w_1}(M(a) - {a_1})}}{{\left\| {M(a) - {a_1}} \right\|}}}\\
	\\
	{ - \frac{{{w_2}(M(a) - {a_2})}}{{\left\| {M(a) - {a_2}} \right\|}}}\\
	\\
	{...}\\
	\\
	{ - \frac{{{w_m}(M(a) - {a_m})}}{{\left\| {M(a) - {a_m}} \right\|}}}
	\end{array}} \right]\begin{array}{*{20}{c}}
{}&{}&{{a_i} \ne M(a),\forall i = 2,...m + 1}
\end{array}\\
\\
\\
\left[ {\begin{array}{*{20}{c}}
	{ - \frac{{{w_1}(M(a) - {a_1})}}{{\left\| {M(a) - {a_1}} \right\|}}}\\
	\\
	{...}\\
	\\
	{...}\\
	{{B_{(0,1)}}}\\
	\\
	{...}\\
	\\
	{ - \frac{{{w_m}(M(a) - {a_m})}}{{\left\| {M(a) - {a_m}} \right\|}}}
	\end{array}} \right]\begin{array}{*{20}{c}}
{}&{}&{{a_j} = 0}
\end{array}
\end{array} \right.
\end{array}\]

\end{Theorem}

{\bf Proof.} Firstly, we have to prove the formula subdifferential of $m$. Since $m$ and $g$ is convex so 
\[\begin{array}{*{20}{l}}
	{{a^*} \in \partial m(a)}\\
	\begin{array}{l}
		\Leftrightarrow \left\langle {{a^*},b - a} \right\rangle  \le m(b) - m(a) = g(M(b),b) - g(M(a),a)\begin{array}{*{20}{c}}
			{}&{\forall b \in \Delta }
		\end{array}\\
		\Leftrightarrow \left\langle {{a^*},b - a} \right\rangle  \le g(y,b) - g(M(a),a)\begin{array}{*{20}{c}}
			{}&{\forall (y,b) \in C \times \Delta }
		\end{array}
	\end{array}\\
	{ \Leftrightarrow \left\langle {(0,{a^*}),(y,b) - (M(a),a)} \right\rangle  \le 0\begin{array}{*{20}{c}}
			{}&{\forall (y,b) \in C \times \Delta }
	\end{array}}\\
	{ \Leftrightarrow (0,{a^*}) \in \partial g(M(a),a)}
\end{array}\]
Next applying Chain rule in Corollary 16.53 of \cite{Bauschke Heinz H. Combettes Patrick L.} and notice that $A$ is surjective, we attain subdifferential of $g$:
\[\begin{array}{ll}
\partial g&= \partial (h \circ A)\\
&= {A^*} \circ \partial h \circ A
\end{array}\]

From Proposition 16.9 of \cite{Bauschke Heinz H. Combettes Patrick L.}, it is possible to calculate subdifferential of $h$ as follows
\[\partial h(y) = \left\{ \begin{array}{l}
{N_C}({y_1}) \times \left\{ {\frac{{{w_1}{y_2}}}{{\left\| {{y_2}} \right\|}}} \right\} \times ... \times \left\{ {\frac{{{w_m}{y_{m + 1}}}}{{\left\| {{y_{m + 1}}} \right\|}}} \right\}\begin{array}{*{20}{c}}
{}&{}&{\text{if }{\forall y_i} \ne 0,i = 2,...m + 1}
\end{array}\\
\\
{N_C}({y_1}) \times ...\left\{ {\frac{{{w_{i - 2}}{y_{i - 1}}}}{{\left\| {{y_{i - 1}}} \right\|}}} \right\} \times \bar B_{(0,1)} \times \left\{ {\frac{{{w_i}{y_{i + 1}}}}{{\left\| {{y_{i + 11}}} \right\|}}} \right\} \times ...\begin{array}{*{20}{c}}
{ \times \left\{ {\frac{{{w_m}{y_{m + 1}}}}{{\left\| {{y_{m + 1}}} \right\|}}} \right\}}&{{\text{if }y_i} = 0}
\end{array}
\end{array} \right.\]

$A$ and $A^*$ may be straightly calculated:
\[A = \left[ {\begin{array}{*{20}{c}}
	{Id}&0&0&0&{...}&0\\
	{Id}&{ - Id}&0&0&{...}&0\\
	{Id}&0&{ - Id}&0&{...}&0\\
	{...}&{...}&{...}&{...}&{...}&{...}\\
	{...}&{...}&{...}&{...}&{...}&{...}\\
	{Id}&{}&{}&{}&{}&{ - Id}
	\end{array}} \right]\]

\[{A^*} = \left[ {\begin{array}{*{20}{c}}
	{Id}&{Id}&{Id}&{Id}&{...}&{Id}\\
	0&{ - Id}&0&0&{...}&0\\
	0&0&{ - Id}&0&{...}&0\\
	{...}&{...}&{...}&{...}&{...}&{...}\\
	{...}&{...}&{...}&{...}&{...}&{...}\\
	0&{}&{}&{}&{}&{ - Id}
	\end{array}} \right]\]

Finally, we obtain subdifferential of $g$ at $(x,a) \in C \times \Delta$
\[\partial g(x,a) = \left\{ \begin{array}{l}
\left[ {\begin{array}{*{20}{c}}
	{{N_C}(x) + \sum\limits_{i = 1}^m {\frac{{{w_i}(x - {a_i})}}{{\left\| {x - {a_i}} \right\|}}} }\\
	{ - \frac{{{w_1}(x - {a_1})}}{{\left\| {x - {a_1}} \right\|}}}\\
	{...}\\
	{ - \frac{{{w_m}(x - {a_m})}}{{\left\| {x - {a_m}} \right\|}}}
	\end{array}} \right]\begin{array}{*{20}{c}}
{}&{}&\text{if }{{x} \ne {a_j},\forall i = 2,...m + 1}
\end{array}\\
\\
\\
\left[ {\begin{array}{*{20}{c}}
	{{N_C}(x) + \sum\limits_{i = 1,i \ne j}^m {\frac{{{w_i}(x - {a_i})}}{{\left\| {x - {a_i}} \right\|}} + {B_{(0,1)}}} }\\
	{ - \frac{{{w_1}(x - {a_1})}}{{\left\| {x - {a_1}} \right\|}}}\\
	{...}\\
	{{B_{(0,1)}}}\\
	{...}\\
	{ - \frac{{{w_m}(x - {a_m})}}{{\left\| {x - {a_m}} \right\|}}}
	\end{array}} \right]\begin{array}{*{20}{c}}
{}&{}&\text{if }{{a_j} = 0}
\end{array}
\end{array} \right.\]
The formular was proven. \qedsymbol
\section{Projected Weiszfeld's algorithm }  
\subsection{Description} \rm
\markboth{\centerline{\sc Introduction}}{\centerline{\sc N.D.Son}} \setcounter{equation}{0}
{\bf Algorithm.} Projected Weiszfeld's algorithm
 \begin{flushleft}
  Initialization: ${x_0} \in C$ and $\epsilon >0$ is tolerance .	\\
  {\bf Step 1: }Compute: 
  \begin{equation} 
  {x_{k + 1}} = {\Pi _C} \circ T({x_k}). 
  \end{equation} 
  {\bf Step 2: }Stop the execution if 
  \begin{equation} 
  \left\| {{x_k} - {x_{k + 1}}} \right\| \le \varepsilon .\label{}
  \end{equation}
  and declare $x_(k+1)$ as solution to the problem (FTW). Otherwise return to Step $1$.
\end{flushleft}

\begin{Remark} \rm \label{remark2.4}
		It is clear that the sequence ${\{ {x_k}\} _{k \ge 0}}$ is a subset of the compact convex set ${\Pi _C}( \overline {co}\{ {a_1},...,{a_m}\} )$. From $(c)$ of Remark \ref{remark1.1} and \eqref{2.4}, then for every $k \ge 0$ we have the result as follows ,
		\begin{equation} \label{2.7}
		\begin{array}{*{20}{c}}
		&{}
		\end{array}\left\langle {T({x_k}) - {x_{k + 1}},x - {x_{k + 1}}} \right\rangle  \le 0\begin{array}{*{20}{c}}
		,&{\text{for all }x \in C.}
		\end{array}
		\end{equation}.
\end{Remark}

\subsection{Convergence to the solution} \rm
We see that proving the convergence of the sequence is not too complicated, but the difficulty is which converging to the unique optimal solution ${\bar x}$ when ${\bar x}$ belongs to set $A$. We will extend the pages of A. Beck and S. Sabach, 2015, \cite{B.A and S.S} for our sequence ${\{ {x_k}\} _{k \ge 0}}$. From now on, we always assume that the sequence $\{ {x_k}\} \notin A$ for every $k \ge 0$.
\medskip 

Consider $C$ is the set in the problem (FTW), we consider the $k$ problem as follows,
\[\min \left\{ {\left. {h(x,{x_k}) = \sum\limits_{i = 1}^m {{w_i}\frac{{{{\left\| {x - {a_i}} \right\|}^2}}}{{\left\| {{x_k} - {a_i}} \right\|}}} } \right|x \in C} \right\},\]
or it can rewrite the $k$ non-constrained problem,
\[\min \left\{ {\left. {h(x,{x_k}) = \sum\limits_{i = 1}^m {{w_i}\frac{{{{\left\| {x - {a_i}} \right\|}^2}}}{{\left\| {{x_k} - {a_i}} \right\|}} + {I_C}(x)} } \right|x \in H} \right\},\]

 \begin{Lemma}\rm \label{lemma3.1}
	The following properties of the $k$ auxiliary functions $h(.,{x_k}):H \to H$ hold:
	\begin{enumerate}[i)]
		\item $h(x_k,x_k)=f(x_k).$
		\item $h(x_{k+1},x_k) \ge 2f(x_{k+1})-f(x_k).$
		\item $x_{k+1} = argmin{_{x \in C}}h(x,x_k).$
	\end{enumerate}
\end{Lemma}
{\bf Proof}.  $i)$ it is trivial.
\medskip

$ii)$ For every two real numbers $a,b>0$, the inequality
\[\frac{{{a^2}}}{b} \ge 2a - b,\]
hold true. Thus, for every $i=1,2,...,m$, $x \in C$ and $x_k \in C \backslash A$, we have
\[\dfrac{{{{\left\| {x_{k+1} - {a_i}} \right\|}^2}}}{{\left\| {x_k - {a_i}} \right\|}} \ge 2\left\| {x_{k+1} - {a_i}} \right\| - \left\| {x_k - {a_i}} \right\|\]
We sum over $i=1,2,...,m$, the result follows.
\medskip

$iii)$ $h(x,x_k)$ is a strict convex function and its unique minimizer is determined by the optimality condiction, 
\[0 \in {\partial _x}h(x,x_k) = 2\sum\limits_{i = 1}^m {\dfrac{{x - {a_i}}}{{\left\| {x_k - {a_i}} \right\|}} + N(x,C),} \]
it is equivalent,
\[0 \in \sum\limits_{i = 1}^m {\dfrac{{x - {a_i}}}{{\left\| {{x_k} - {a_i}} \right\|}}}  + N(x,C),\]
or
\[T({x_k}) \in x + N(x,C),\]
so
\[x = \Pi  \circ T({x_k}),\]
it can be seen that for any $k \ge 0$,
\[x_{k+1} = argmin_{{x \in C}}h(x,x_k).\]

We are now able to prove the  descent property of the algorithm.
\begin{Lemma} \label{lemma3.2}
For every $k \ge 0$ such that $x_k \notin A$, we have
\begin{equation} \label{3.7}
f({x_{k + 1}}) \le f({x_k}),
\end{equation}
moreover,
\begin{equation} \label{3.8}
	f({x_{k + 1}}) = f({x_k}) \text{ if and only if } {x_{k + 1}} = {x_k}.
	\end{equation}
\end{Lemma} 
{\bf Proof}. From $iii)$ of Lemma \ref{lemma3.1} and the strict convexity of the $k$ auxiliary functions $h(x,x_k)$, we have
\begin{equation} \label{3.9}
h({x_{k + 1}},{x_k}) < h({x_k},{x_k}) = f({x_k}).
\end{equation}
From $ii)$ of Lemma \ref{lemma3.1}, we obtain
\begin{equation} \label{3.10}
h({x_{k + 1}},{x_k}) \ge 2f({x_{k + 1}}) - f({x_k}).
\end{equation}
From the formulas \eqref{3.9} and \eqref{3.10}, it implies \eqref{3.7}. And if $f({x_{k + 1}}) = f({x_k})$ and $x_{k + 1}  \ne  {x_k}$ , then we easily have a contradiction.
		
\begin{Lemma} \label{lemma3.3}
	Suppose that $x_k \notin A$ for all $k \ge 0$. Let $L({x_k}) = \sum\limits_{i = 1}^m {\dfrac{{{w_i}}}{{\left\| {{x_k} - {a_i}} \right\|}}}$, then
	\begin{equation} \label{3.11}
	f({x_{k + 1}}) \le f({x_k}) + \dfrac{{L({x_k})}}{2}\left[ {{{\left\| {{x_{k + 1}} - {x_k}} \right\|}^2} + 2\left\langle {{x_k} - T({x_k}),{x_{k + 1}} - {x_k}} \right\rangle } \right].
	\end{equation}	
\end{Lemma}

{\bf Proof}. For every $i=1,2,...,m$, we have
\[{w_i}\dfrac{{{{\left\| {x_{k+1} - {a_i}} \right\|}^2}}}{{\left\| {{x_k} - {a_i}} \right\|}} = {w_i}\dfrac{{{{\left\| {x_{k+1} - {x_k} + {x_k} - {a_i}} \right\|}^2}}}{{\left\| {{x_k} - {a_i}} \right\|}}\]
Use the property of $<.,.>$ in the real Hilbert space $H$, it implies
\[\begin{array}{*{20}{l}}
{{w_i}\dfrac{{{{\left\| {{x_{k + 1}} - {x_k} + {x_k} - {a_i}} \right\|}^2}}}{{\left\| {{x_k} - {a_i}} \right\|}}}&{ = {w_i}\dfrac{{{{\left\| {{x_{k + 1}} - {x_k}} \right\|}^2} + {{\left\| {{x_k} - {a_i}} \right\|}^2} + 2\left\langle {{x_k} - {a_i},{x_{k + 1}} - {x_k}} \right\rangle }}{{\left\| {{x_k} - {a_i}} \right\|}}}\\
{}&{ = {w_i}\dfrac{{{{\left\| {{x_{k + 1}} - {x_k}} \right\|}^2}}}{{\left\| {{x_k} - {a_i}} \right\|}} + 2\left\langle {{w_i}\dfrac{{{x_k} - {a_i}}}{{\left\| {{x_k} - {a_i}} \right\|}},{x_{k + 1}} - {x_k}} \right\rangle  + {w_i}\left\| {{x_k} - {a_i}} \right\|.}
\end{array}\]
Summing over $i=1,2,...,m,$ the result follows
\[\begin{array}{*{20}{l}}
{h({x_{k + 1}},{x_k})}&{ = f({x_k}) + \sum\limits_{i = 1}^m {{w_i}\dfrac{{{{\left\| {{x_{k + 1}} - {x_k}} \right\|}^2}}}{{\left\| {{x_k} - {a_i}} \right\|}}}  + 2\left\langle {\sum\limits_{i = 1}^m {{w_i}\dfrac{{{x_k} - {a_i}}}{{\left\| {{x_k} - {a_i}} \right\|}}} ,{x_{k + 1}} - {x_k}} \right\rangle ,}\\
{}&{}\\
{}&{ = f({x_k}) + L({x_k})\left[ {{{\left\| {{x_{k + 1}} - {x_k}} \right\|}^2} + 2\left\langle {{x_k} - T({x_k}),{x_{k + 1}} - {x_k}} \right\rangle } \right].}
\end{array}\]
Applying ii) of Lemma \ref{lemma3.1} yelds
\[\begin{array}{*{20}{l}}
{\begin{array}{*{20}{c}}
	{}&{2f({x_{k + 1}}) - f({x_k}) \le h({x_{k + 1}},{x_k}),}
	\end{array}}\\
{}\\
{ \Leftrightarrow 2f({x_{k + 1}}) \le 2f({x_k}) + L({x_k})\left[ {{{\left\| {{x_{k + 1}} - {x_k}} \right\|}^2} + 2\left\langle {{x_k} - T({x_k}),{x_{k + 1}} - {x_k}} \right\rangle } \right],}\\
{}\\
{ \Leftrightarrow f({x_{k + 1}}) \le f({x_k}) + \dfrac{{L({x_k})}}{2}\left[ {{{\left\| {{x_{k + 1}} - {x_k}} \right\|}^2} + 2\left\langle {{x_k} - T({x_k}),{x_{k + 1}} - {x_k}} \right\rangle } \right].}
\end{array}\]

\begin{Lemma}\label{lemma3.4}
	For any $x \in C$, we have the following inequality holds
	\begin{equation} \label{3.12}
	f({x_{k + 1}}) - f(x) \le \dfrac{{L({x_k})}}{2}\left( {{{\left\| {{x_k} - x} \right\|}^2} - {{\left\| {{x_{k + 1}} - x} \right\|}^2}} \right).
	\end{equation}
\end{Lemma}
{\bf Proof}. We consider the function $g(y) = \sum\limits_{i = 1}^m {{w_i}\left\| {{y} - {a_i}} \right\|}$ with for any $y \in H$. By the gradient inequality of the function $g(.)$, for any $x \in C$ we have
\[g({x_k}) \le g(x) + \left\langle {\sum\limits_{i = 1}^m {{w_i}\dfrac{{{x_k} - {a_i}}}{{\left\| {{x_k} - {a_i}} \right\|}}} ,{x_k} - x} \right\rangle, \]
or
\[g({x_k}) - g(x) \le L({x_k})\left\langle {{x_k} - T({x_k}),{x_k} - x} \right\rangle \]

Since $g(x)=f(x)$ for every $x \in C$ and \eqref{3.11} of Lemma \ref{lemma3.3}, we obtain
\[\begin{array}{*{20}{l}}
{f({x_{k + 1}}) + f({x_k}) - f(x) \le f({x_k}) + \dfrac{{L({x_k})}}{2}}&{\left[ {{{\left\| {{x_{k + 1}} - {x_k}} \right\|}^2} + 2\left\langle {{x_k} - T({x_k}),{x_{k + 1}} - {x_k}} \right\rangle } \right.}\\
{}&{\left. { + 2\left\langle {{x_k} - T({x_k}),{x_k} - x} \right\rangle } \right],}
\end{array}\]
 so
 \begin{equation}\label{3.13}
f({x_{k + 1}}) - f(x) \le \dfrac{{L({x_k})}}{2}\left[ {{{\left\| {{x_{k + 1}} - {x_k}} \right\|}^2} + 2\left\langle {{x_k} - T({x_k}),{x_{k + 1}} - x} \right\rangle } \right].
 \end{equation}
 
 Use the property of $<.,.>$ as follows,
 \begin{eqnarray} \label{3.14}
 \begin{array}{ll}
 {\left\| {{x_{k + 1}} - {x_k}} \right\|^2} &= {\left\| {{x_{k + 1}} - x + x - {x_k}} \right\|^2}\\
 &= {\left\| {{x_{k + 1}} - x} \right\|^2} + {\left\| {{x_k} - x} \right\|^2} + 2\left\langle {x- {x_k},{x_{k + 1}} - x} \right\rangle .
 \end{array}
 \end{eqnarray}
 From the formulation (\ref{3.13}) and (\ref{3.14}), we obtain
 \[\begin{array}{ll}
 f({x_{k + 1}}) - f(x) &\le \dfrac{{L({x_k})}}{2}\left[ {{{\left\| {{x_{k + 1}} - x} \right\|}^2} + {{\left\| {{x_k} - x} \right\|}^2} + 2\left\langle {{x_k} - T({x_k}) + x - {x_k},{x_{k + 1}} - x} \right\rangle } \right],\\
 \\
 &\le \dfrac{{L({x_k})}}{2}\left[ {{{\left\| {{x_{k + 1}} - x} \right\|}^2} + {{\left\| {{x_k} - x} \right\|}^2} + 2\left\langle {x - T({x_k}),{x_{k + 1}} - x} \right\rangle } \right],\\
 \\
 &\le \dfrac{{L({x_k})}}{2}\left[ {{{\left\| {{x_{k + 1}} - x} \right\|}^2} + {{\left\| {{x_k} - x} \right\|}^2} + 2\left\langle {x - {x_{k + 1}} + {x_{k + 1}} - T({x_k}),{x_{k + 1}} - x} \right\rangle } \right],\\
 \\
 &\le \dfrac{{L({x_k})}}{2}{\left[ {\left\| {{x_{k + 1}} - x} \right\| - \left\| {{x_k} - x} \right\|} \right]^2} + L({x_k})\left\langle {T({x_k}) - {x_{k + 1}},x - {x_{k + 1}}} \right\rangle .
 \end{array}\]
 Repeat Remark \ref{remark2.2},
 \[\left\langle {T({x_k}) - {x_{k + 1}},x - {x_{k + 1}}} \right\rangle  \le 0\begin{array}{*{20}{c}}
 {}&{\text{for any $x \in C$}}
 \end{array},\]
 so the proof is completed.
 
 \begin{Corollary} \label{}
 	If we have a point $x \in C$ such that $f(x) \le f({x_k})\begin{array}{*{20}{c}}
 	{\text{for all}}&{k \ge 0,}
 	\end{array}$ then the left-side of (\ref{3.12}) is nonnegative. We have the following inequation,
 	\begin{equation}
 	\left\| {{x_{k + 1}} - x} \right\| \le \left\| {{x_k} - x} \right\|.
 	\end{equation} \label{}
 	Specially, we replace $x$ by the optimal solution ${\bar x}$,
 	\begin{equation}
 	\left\| {{x_{k + 1}} - {\bar x}} \right\| \le \left\| {{x_k} - {\bar x}} \right\|.
 	\end{equation} \label{}
 	Suppose existence $l \ge 0$ such that $\left\| {{x_{l + 1}} - {{\bar x}}} \right\| = \left\| {{x_l} - {{\bar x}}} \right\|$, then from the \eqref{3.8} of Lemma \ref{lemma3.2}, we get $x_l={\bar x}$.
 \end{Corollary}
 
From the above-lemma, we prove that the sequence ${\{ {x_k}\} _{k \ge 0}}$ converges to the unique solution of the problem (FTW).
\begin{Theorem}
	Suppose that for any $k \ge 0$, $x_k \notin A$ then the sequence ${\{ {x_k}\}}$ converges to the unique solution ${{\bar x} }$of (FTW) .
\end{Theorem}
 {\bf Proof.} First, we will prove the sequence ${\{ {x_k}\} _{k \ge 0}}$ is the convergence sequence. Indeed, repeat Remark \ref{remark2.4}, the sequence ${\{ {x_k}\} _{k \ge 0}}$ is a subset of the compact set $ \Pi_C (\Omega )$, so we can take two subsequence ${\{ {x_{{k_l}}}\} _{l \ge 0}}$ and ${\{ {x_{{k_m}}}\} _{m \ge 0}}$ converging to limits $\bar x$ and $\tilde x$, respectively. From Lemma \ref{lemma3.4}, it follows that $f(\bar x) \le f({x_k})$ for all $k \ge 0$, and thus from Corollary \ref{}, we get that the sequence ${\{ \left\| {{x_k} - \bar x} \right\|\} _{k \ge 0}}$ is nonincreasing, it thus coverges to some scalar $r$. It is clear that
 \[r = \mathop {\lim }\limits_{k \to \infty } \left\| {{x_k} - \bar x} \right\| = \mathop {\lim }\limits_{l \to \infty } \left\| {{x_l} - \bar x} \right\| = 0,\]\
 but, on the other hand 
 \[r = \mathop {\lim }\limits_{k \to \infty } \left\| {{x_k} - \bar x} \right\| = \mathop {\lim }\limits_{m \to \infty } \left\| {{x_m} - \bar x} \right\| = \left\| {\tilde x - \bar x} \right\|.\]
So $\left\| {\tilde x - \bar x} \right\| = 0$, which leading to that the sequence  ${\{ {x_k}\}}_{k \ge 0}$ converges. Second, we assume that the sequence ${\{ {x_{{k}}}\} _{k \ge 0}}$ converges to some ${\bar x}$ then ${\bar x}$ is the optimal solution $(FTW)$ of the problem (FTW). Indeed, if ${\bar x}$ is not in the set $A$, since ${x_{k + 1}} = {\Pi _C} \circ T({x_k})$ and the mapping $\Pi_C  \circ T(.)$ is continuous, so ${\bar x} = {\Pi _C} \circ T({{\bar x}})$ . From Theorem \ref{theorem2.1} and the formula \eqref{2.3} then ${\bar x}$ is $x ^ *$. If ${\bar x}$ is in the set $A$, supposing that ${\bar x} = a_j$ with for some $j \in \{ 1,2,...,m\} $. From $iii)$ of Lemma \ref{lemma3.1}, we have
  \[0 \in {\partial _x}h({x_{k + 1}},{x_k}) = 2\sum\limits_{i = 1}^m {\frac{{{x_{k + 1}} - {a_i}}}{{\left\| {{x_k} - {a_i}} \right\|}}}  + N({x_{k + 1}},C),\]
  or
\begin{equation}\label{3.17}
{ - \sum\limits_{i = 1}^m {\frac{{{x_{k + 1}} - {a_i}}}{{\left\| {{x_k} - {a_i}} \right\|}}}  \in N({x_{k + 1}},C).}
\end{equation}
We fix an abarary element $x$ which belongs to $C$ set. Then, \eqref{3.17} is equivalent
\[\begin{array}{ll}
&\left\langle { - \sum\limits_{i = 1}^m {\dfrac{{{x_{k + 1}} - {a_i}}}{{\left\| {{x_k} - {a_i}} \right\|}}} ,x - {x_{k + 1}}} \right\rangle  \le 0,\\
\\
\Leftrightarrow & \left\langle { - \dfrac{{{x_{k + 1}} - {a_j}}}{{\left\| {{x_k} - {a_j}} \right\|}} - \sum\limits_{i = 1,i \ne j}^m {\dfrac{{{x_{k + 1}} - {a_i}}}{{\left\| {{x_k} - {a_i}} \right\|}}} ,x - {x_{k + 1}}} \right\rangle  \le 0,\\
\end{array}\]
so we get the following formula,
\begin{equation} \label{3.18}
\left\langle { - \dfrac{{{x_{k + 1}} - {a_j}}}{{\left\| {{x_k} - {a_j}} \right\|}},x - {x_{k + 1}}} \right\rangle  + \left\langle { - \sum\limits_{i = 1,i \ne j}^m {\dfrac{{{x_{k + 1}} - {a_i}}}{{\left\| {{x_k} - {a_i}} \right\|}}} ,x - {x_{k + 1}}} \right\rangle  \le 0,
\end{equation}
From Corollary \ref{}, we see that the sequence ${\left\{ {\dfrac{{{x_{k + 1}} - {a_j}}}{{{x_k} - {a_j}}}} \right\}_{k \ge 0}}$ is subset of the unit circle ${{\bar B}_{({0},1)}}$. Since $H$ is a real Hilbert space, so it exists a subsequence ${\left\{ {\dfrac{{{x_{{k_l} + 1}} - {a_j}}}{{{x_{{k_l}}} - {a_j}}}} \right\}_{l \ge 0}}$ converges weakly to some $u \in {{\bar B}_{({0_H},1)}}$ (convergence in weak-topology). Since the sequence $\{x_k\}$ converge to $a_j$ (convergence in the norm), so
\[\left\langle { - \dfrac{{{x_{{k_l} + 1}} - {a_j}}}{{\left\| {{x_{{k_l}}} - {a_j}} \right\|}},x - {x_{{k_l} + 1}}} \right\rangle  \to \left\langle { - u,x - {a_j}} \right\rangle \text{ when $l \to \infty$}.\]
It is easy that 
\[\left\langle { - \sum\limits_{i = 1,i \ne j}^m {\dfrac{{{x_{{k_l} + 1}} - {a_i}}}{{\left\| {{x_{{k_l}}} - {a_i}} \right\|}}} ,x - {x_{{k_l} + 1}}} \right\rangle  \to \left\langle { - \sum\limits_{i = 1,i \ne j}^m {\dfrac{{{a_j} - {a_i}}}{{\left\| {{a_j} - {a_i}} \right\|}}} ,x - {a_j}} \right\rangle \text{ when $l \to \infty$}.\]
Hence when $l \to \infty$, \eqref{3.18} becomes
\[\left\langle { - \sum\limits_{i = 1,i \ne j}^m {\dfrac{{{a_j} - {a_i}}}{{\left\| {{a_j} - {a_i}} \right\|}}}  - u,x - {a_j}} \right\rangle  \le 0.\]
From Theorem \ref{theorem2.1}, we show that $a_j$ is an optimal solution of (FTW). Proof is completed. \qedsymbol

We can ignore the part proof of the algorithm's convergence to $a_j$ when it is a solution by the following way. Firstly, we check optimal condition Theorem 2.2 \ref{theorem2.1} ii) at $a1,...,a_m$ and if no $a_i$ satisfies, we apply the algorithm. However, checking the optimal condition at $a1,...,a_m$ is not easy unless additionally $C$ is either cone or $int(C) \ne empty$ as seeing in Theorem 2.2 \ref{theorem2.1} ii).

\end{document}